\documentclass[a4paper,10pt]{article}
\usepackage{amsfonts}
\usepackage{bbm}
\usepackage{mathrsfs}
\usepackage{amsmath,amsthm,amssymb,amscd}
\usepackage[center]{titlesec}
\newtheorem{theo}{Theorem}[section]
\newtheorem{prop}[theo]{Proposition}
\newtheorem{lem}[theo]{Lemma}
\newtheorem{rem}[theo]{Remark}

\newtheorem{defi}[theo]{Definition}

\newcommand{\bp}{\begin{proof}}
\newcommand{\ep}{\end{proof}}

\begin{document}
 \setlength{\baselineskip}{13pt} \pagestyle{myheadings}

 \title{
 {On cohomology of almost complex 4-manifolds}
 \thanks{
 Supported by NSFC (China) Grants 11071208 (Wang), 10871139, 11271276 (Zhang) and 11101352 (Zhu).\, %
 Correspondence to: hywang@yzu.edu.cn (Wang).}
 }
 \author{{Qiang Tan, Hongyu Wang, Ying Zhang, and Peng Zhu}
 }
 \date{}
 \maketitle

 \noindent {\bf Abstract}
 Based on recent work of T. Draghici, T.-J. Li and W. Zhang,
 we further investigate properties of the dimension $h_J^-$ of the $J$-anti-invariant cohomology
 subgroup $H_J^-$ of a closed almost Hermitian 4-manifold $(M,g,J,F)$
 using metric compatible almost complex structures. 
 We prove that $h_J^-=0$ for generic almost complex structures $J$ on $M$.
 \\

\noindent {{\bf Keywords} Almost Hermitian 4-manifold $\cdot$ $J$-anti-invariant cohomology $\cdot$
            Metric compatible almost complex structure}\\
            

 \section{Introduction}

 For an almost complex manifold $(M,J)$,
 T.-J. Li and W. Zhang \cite{LiZh} introduced subgroups $H_{J}^{+}$ and $H_{J}^{-}$ %
 of the real degree 2 de Rham cohomology group $H^2(M,\mathbb R)$, as the sets of cohomology classes %
 which can be represented by $J$-invariant and $J$-anti-invariant real 2-forms, respectively. %
 Let us denote by $h_J^{+}$ and $h_J^{-}$ the dimensions of $H_{J}^{+}$ and $H_{J}^{-}$, respectively.

 It is interesting to consider whether or not the subgroups $H_{J}^{+}$ and $H_{J}^{-}$ %
 induce a direct sum decomposition of $H^2(M,\mathbb R)$.
 In the case of direct sum decomposition, $J$ is said to be $C^\infty$ pure and full (see Definition \ref{2000}).

 This is known to be true for integrable almost complex structures $J$
 which admit compatible K\"{a}hler metrics on compact manifolds of any dimension.
 In this case, the induced decomposition is nothing but the classical real
 Hodge-Dolbeault decomposition of $H^2(M,\mathbb R)$ (see \cite{Barth}).

 Note that there are topological obstructions to the existence of
 almost complex structures on an even dimensional manifold.
 For a closed 4-manifold, a necessary condition is that $1-b_1+b^+$ be even \cite{Barth}, %
 where $b_1$ is the first Betti number and $b^+$ is the number of positive eigenvalues of
 the quadratic form on $H^2(M,\mathbb R)$ defined by the cup product, %
 hence the condition is either $b_1$ be even and $b^+$ odd, or $b_1$ be odd and $b^+$ even. %

 It is a well-known fact that any closed complex surface with $b^+$ odd is K\"{a}hler. %
 This was originally obtained from the classification theory, but
 direct proofs have been given in \cite{Bu,La}.

 In dimension 4, T. Draghici, T.-J. Li and W. Zhang \cite{DLZ1}
 proved that on any closed almost complex 4-manifold $(M,J)$, $J$ is $C^\infty$ pure and full.
 Further in \cite{DLZ2}, they computed the subgroups $H_{J}^{+}$ and $H_{J}^{-}$ %
 and their dimensions $h_J^{+}$ and $h_J^{-}$ for almost complex structures metric related to an integrable one.
 Using Gauduchon metrics (\cite{G2}), they proved that the almost complex structures $\tilde{J}$ with $h_{\tilde{J}}^{-}=0$ %
 form an open dense set in the $C^{\infty}$-Fr\'{e}chet-topology in the space
 of almost complex structures metric related to an integrable one (\cite[Theorem 1.1]{DLZ2}). %
 Based on this, they made a conjecture (Conjecture 2.4 in \cite{DLZ2}) about the dimension $h_J^-$ of
 $H_J^-$ on a compact 4-manifold which asserts that $h_J^-$ vanishes for
 generic almost complex structures $J$.
 In particular, they have confirmed their conjecture for 4-manifolds with $b^+=1$
 (\cite[Theorem 3.1]{DLZ2}).

 In this paper we confirm the conjecture completely (see Theorem \ref{thm1} below)
 by using $g$-compatible almost complex structures.
%

 Let $\mathcal{J}$ be the space of all almost complex structures on $M$ %
 and denote by $\mathcal{J}_{g}^{\rm c}$ and $\mathcal{J}_{F}^{\rm t}$ respectively %
 the spaces of $g$-compatible and $F$-tame almost complex structures on $M$; namely, %
 \begin{align*}
 & \mathcal{J}_{g}^{\rm c} = \{ J\in\mathcal{J} \mid g(JX,JY)=g(X,Y), \forall X,Y \in TM \}, \\ %
 & \mathcal{J}_{F}^{\rm t} = \{ J\in\mathcal{J} \mid F(X,JX)>0, \forall X \in TM, X \neq 0 \}.
 \end{align*}
 It is well known that $\mathcal{J}_{g}^{\rm c}$ and $\mathcal{J}_{F}^{\rm t}$ are
 contractible $C^\infty$-Fr\'{e}chet spaces.
 See \cite{Au, DLZ1, DLZ2} for details. 

%
%

 \begin{theo}\label{thm1}
 Let $M$ be a closed 4-manifold admitting almost complex structures.
 Then the set of almost complex structures $J$ on $M$ with $h_J^-=0$
 is an open dense subset of $\mathcal{J}$ in the $C^\infty$-topology.
 \end{theo}


 The rest of the paper is organized as follows.
 In \S 2 we recall definitions, preliminary results mainly as given in \cite{DLZ2},
 and constructions of $g$- and $F$-compatible almost complex structures.
 Finally in \S 3 we give the proof of Theorem \ref{thm1}.


 \section{Definitions and Preliminaries}\setcounter{equation}{0}

 Suppose $(M,J)$ is a closed almost complex 4-manifold.
 One can construct a $J$-invariant Riemannian metric $g$ on $M$.
 Such a metric $g$ is called an almost Hermitian metric for $(M,J)$.
 This then in turn gives a $J$-compatible non-degenerate 2-form $F$ by $F(X,Y)=g(JX,Y)$, %
 called the fundamental 2-form.
 Such a quadruple $(M,g,J,F)$ is called a closed almost Hermitian 4-manifold. %
 Thus an almost Hermitian structure on $M$ is a triple $(g,J,F)$. %
 By direct calculation, $F \wedge F = 2d\mu_g$,
 where $d\mu_g$ is the volume form of $M$ determined by $g$. %


 Note that $J$ acts on the space $\Omega^2$ of 2-forms on $M$ as an involution by
 \begin{align}
 \alpha \longmapsto \alpha(J\cdot,J\cdot), \quad \alpha\in\Omega^2(M).
 \end{align}
%
 This gives the $J$-invariant, $J$-anti-invariant decomposition of 2-forms:
 \begin{align}
 \Omega^2 = \Omega^+_J \oplus \Omega^-_J, \quad \alpha = \alpha_J^+ + \alpha_J^- %
 \end{align}
 as well as the splitting of corresponding vector bundles
 \begin{align}\label{201}
 {\Lambda}^2={\Lambda}_J^+ \oplus {\Lambda}_J^-. %
 \end{align}

 \begin{rem}
 {\rm
 Note that $\Lambda_J^-$ is a vector bundle of rank two and
 $\Lambda_J^-$ inherits an almost complex structure, still denoted by $J$,
 defined by
 \begin{align*}
 (J\alpha)(\cdot,\cdot) = -\alpha(J\cdot,\cdot), \quad\quad \alpha \in \Lambda_J^-.
 \end{align*}
 It is well known that, when $J$ is integrable,
 $\beta\in \mathcal Z_J^-$ if and only if $J\beta\in \mathcal Z_J^-$.
 Conversely, if $(M,J)$ is a connected almost complex 4-manifold and
 there exists nonzero $\beta\in \mathcal Z_J^-$ such that $J\beta\in \mathcal Z_J^-$,
 then $J$ is integrable (see \cite{Sa}).
 }
 \end{rem}


 \begin{defi}\label{2000}{\rm (cf. \cite{DLZ1,LiZh}) %
 Let $\mathcal Z^2$ denote the space of closed 2-forms on $M$ and set
 \begin{align*}
 \mathcal Z_J^+ := \mathcal Z^2 \cap \Omega_J^+, \quad
 \mathcal Z_J^- := \mathcal Z^2 \cap \Omega_J^-.
 \end{align*}
 Define the $J$-invariant and $J$-anti-invariant cohomology subgroups $H_J^{\pm}$ by %
 \begin{align*}
 H^\pm_J=\{\mathfrak a \in H^2(M;\mathbb{R}) \mid %
           \mbox{there exists} \ \alpha \in \mathcal Z_J^{\pm} \ \mbox{such that} \ \mathfrak a = [\alpha]\}. %
 \end{align*}
 We say $J$ is $C^\infty$ {\it pure} if $H_J^ + \cap H_J^- = \{0\}$, $C^\infty$ {\it full} if $H_J^+ + H_J^- = H^2(M;\mathbb{R})$, %
 and $J$ is $C^\infty$ {\it pure and full} if %
 \begin{align*}
  H^2(M;\mathbb{R}) = H_J^+\oplus H_J^-.
 \end{align*} }
 \end{defi}

 \begin{prop} {\rm (cf. \cite[Theorem 2.3]{DLZ1})}\label{prop:2.3}
 For any closed almost complex 4-manifold $(M,J)$, $J$ is $C^\infty$ pure and full.
 \end{prop}


 Since $(M,g,J,F)$ is a closed almost Hermitian 4-manifold, the Hodge star operator $*_g$ gives
 the well-known self-dual, anti-self-dual decomposition of 2-forms as well as the corresponding splitting of the bundle (see \cite{DK}):
 \begin{equation}\label{202a}
 \Omega^2 = \Omega_g^+ \oplus \Omega_g^-, \quad \alpha = \alpha_g^+ + \alpha_g^-;
 \end{equation}
 \begin{equation}\label{202}
 {\Lambda}^2 = {\Lambda}_g^+ \oplus {\Lambda}_g^-.
 \end{equation}
%
 Since the Hodge-de Rham-Laplace operator commutes with $*_g$,
 the decomposition \eqref{202} holds for the space $\mathcal {H}_g$ of harmonic 2-forms as well.
 By Hodge theory, this induces cohomology decomposition by the metric $g$:
 \begin{align}\label{203}
 H^2(M;\mathbb{R})=\mathcal {H}_g=\mathcal {H}_g^+\oplus\mathcal
 {H}_g^-.
 \end{align}
 Similar to Definition \ref{2000}, one defines
 \begin{align}
 H_g^{\pm} = \{ \mathfrak a \in H^2(M;\mathbb{R}) \mid
             \mathfrak a = [\alpha] \,\; \mbox{for some} \,\;\alpha \in \mathcal Z_g^{\pm}:=\mathcal Z^2\cap\Omega_g^{\pm} \}. %
 \end{align}
%
%
%
 It is easy to see that %
 \begin{align*}
 H_g^\pm = \mathcal Z_g^{\pm} = \mathcal {H}_g^\pm
 \end{align*}
 and \eqref{203} can be written as
 \begin{align}
 H^2(M;\mathbb{R}) = H_g^+ \oplus H_g^-.
 \end{align}
 There are the following relations between the decompositions \eqref{201} and \eqref{202}
 on an almost Hermitian 4-manifold:
 \begin{align}
 & {\Lambda}_J^+=\mathbb{R} F\oplus{\Lambda}_g^-, \label{205} \\
 & {\Lambda}_g^+=\mathbb{R} F\oplus{\Lambda}_J^-, \label{206} \\
 & {\Lambda}_J^+\cap{\Lambda}_g^+=\mathbb{R} F, \ \ \ \ {\Lambda}_J^-\cap{\Lambda}_g^-=\{0\}. \label{207}
 \end{align}
 See \cite{D3} for details. It is easy to see that
 $H_J^-\subset \mathcal {H}_g^+$ and $\mathcal {H}_g^-\subset H_J^+$.

 \vskip 6pt

 Let $b_2$, $b^{+}$ and $b^{-}$ be the second, the self-dual and
 the anti-self-dual Betti number of $M$, respectively.
 Thus $b_2=b^{+}+b^{-}$. %
%
 It is easy to see that, for a closed almost Hermitian 4-manifold $(M,g,J,F)$, there hold (see \cite{DLZ1}):
 \begin{align}
 H_J^- = \mathcal Z_J^-, \quad
 h_J^+ + h_J^- = b_2, \quad
 h_J^+ \geq b^-, \quad h_J^- \leq b^+.
 \end{align}

 Lejmi recognizes $\mathcal Z_J^-$ as the kernel of an elliptic operator on $\Omega_J^-$.

 \begin{lem}\label{31}{\rm (cf. \cite{Le1,Le2})}
 Let $(M,g,J,F)$ be a closed almost Hermitian 4-manifold.
 Let operator $P: \Omega_J^-\rightarrow\Omega_J^-\nonumber$ be defined by
 \begin{align*}
  P(\psi) = P_J^-(d\delta_g\psi), 
 \end{align*}
 where $P_J^-: \Omega^2 \rightarrow \Omega_J^-$ is the projection, $\delta_g$ is the codifferential operator with respect to metric $g$.
 Then $P$ is a self-adjoint strongly elliptic linear operator with kernel the $g$-harmonic $J$-anti-invariant 2-forms.
 \end{lem}

 Hence one has the decomposition of $\Omega_J^-$:
 \begin{align*}
 \Omega_J^- = {\rm ker}\, P \oplus P_J^-(d\Omega^1) = H_J^- \oplus P_J^-(d\Omega^1). %
 \end{align*}

%
%
%

 Let $H_J^{-,\perp}$ denote the subgroup of $\mathcal {H}_g^+$ which
 is orthogonal to $H_J^-$ with respect to the cup product; that is, %
 \begin{align}
 H_J^{-,\perp} := \{ \omega \in \mathcal Z_g^+ \mid
    \textstyle\int_M \omega\wedge\alpha=0 \ \ \forall \alpha \in \mathcal Z_J^- \}.
 \end{align}

 By Lemma \ref{31} and the results in \cite[Lemmas 2.4 and 2.6]{DLZ1}, we have
 the following lemma which will be used in \S 3.

 \begin{lem}\label{32}
 Let $(M,g)$ be a closed Riemannian 4-manifold.
 If $\alpha \in \Omega_g^+$ and
 $\alpha = \alpha_{\rm h} + d\theta + \delta_g\psi$ is its Hodge decomposition, %
 then $P_g^+(d\theta) =  P_g^+(\delta_g\psi)$ and $P_g^-(d\theta) = -P_g^-(\delta_g\psi)$, %
 where $P_g^{\pm}: \Omega^2 \rightarrow \Omega_g^{\pm}$ are the projections.
 Moreover, the 2-form $\alpha-2P_g^+(d\theta)=\alpha_{\rm h}$ is harmonic and
 $\alpha + 2P_g^-(d\theta) = \alpha_{\rm h} + 2d\theta$.
 In particular, if $(M,g,J,F)$ is a closed almost Hermitian 4-manifold and if $\alpha \in H_J^{-,\perp}$
 is a self-dual harmonic 2-form, then $\alpha = fF + P_J^-(d\theta)$
 for some function $f \not\equiv 0$ and $\alpha-d\theta\in \mathcal Z_J^+$. %
 \end{lem}

 \begin{rem}\label{34}{\rm
 As direct consequences of Lemmas \ref{31} and \ref{32}, we have decompositions
 as self-dual harmonic 2-forms and as cohomology classes: %
 \begin{align*}
 \mathcal {H}_g^+ = H_J^- \oplus H_J^{-,\perp}, \quad\quad %
          H_J^{+} = H_J^{-,\perp} \oplus \mathcal {H}_g^-. %
 \end{align*} }
 \end{rem}

%
%

%

 The almost complex structure $J$ acts on the cotangent bundle %
 $T^*M$ by $(J\alpha)(\cdot) = -\alpha(J\cdot)$, where $\alpha$ is a 1-form.
 Hence $J$ induces an action $J \otimes J$ on $\otimes^2 T^*M$, still denoted by $J$. %
 The bundle of 2-forms decomposes under the action of $J$ as
 \begin{align}
 \Lambda^2 = \Lambda_J^+ \oplus \Lambda_J^-, %
 \end{align}
 and the symmetric tensor bundle of type $(2,0)$ decomposes as
 \begin{align}
 S^2 = S_J^+ \oplus S_J^-. %
 \end{align}

 \vskip 3pt

 Let us construct $g$-compatible almost complex structures
 using sections of the bundle of $J$-anti-invariant 2-forms. %
 Given $\alpha \in \Omega_J^-$, define tensor field $K_{\alpha}$ of type $(1,1)$ by %
 \begin{align}
 g(X,K_{\alpha}Y) = \alpha(X,Y).
 \end{align}
 It can be checked that $K_{\alpha}$ and $JK_{\alpha}$ are skew-adjoint.
 It follows that ${\rm Id}+JK_{\alpha}$ is invertible (\cite{Lee}).
 Define tensor field $J_{\alpha}$ of type $(1,1)$ and 2-form $F_{\alpha}$ by
 \begin{align}\label{Jalpha}
 J_{\alpha} := ({\rm Id}+JK_{\alpha})^{-1} J ({\rm Id}+JK_{\alpha}), \quad
 F_{\alpha} := g(J_{\alpha}\cdot,\cdot).
 \end{align}
 It is easy to see that $J_\alpha^2=-{\rm Id}$;
 hence $J_\alpha$ is an almost complex structure. %

 The $J_{\alpha}$ and $F_{\alpha}$ so defined have the following properties.

 \begin{prop}\label{prop3.1}{\rm (\cite[Proposition 1.5]{Lee})}
 Given $\alpha \in \Omega_J^-$, we define the norm function $|\alpha| \in C^{\infty}(M)$ %
 of $\alpha$ by $\alpha\wedge\alpha = 2 |\alpha|^2 d\mu_g$, i.e., %
 \begin{align*}
 |\alpha|^2 = (\alpha\wedge\alpha)/(2 d\mu_g). %
 \end{align*}
 Then the $J_{\alpha}$ and $F_{\alpha}$ defined by (\ref{Jalpha}) satisfy: %
 \begin{align*}
 & g(J_{\alpha}X,J_{\alpha}Y)=g(X,Y), \\
 & J_{\alpha} = \frac{1-|\alpha|^2}{1+|\alpha|^2} J - \frac{2}{1+|\alpha|^2} K_{\alpha}, \\ %
 & F_{\alpha} = \frac{1-|\alpha|^2}{1+|\alpha|^2} F + \frac{2}{1+|\alpha|^2} \alpha. %
 \end{align*} %
 \end{prop}

 Thus the triple $(g, J_\alpha, F_\alpha)$ is an almost Hermitian structure on $M$.

 \begin{rem}{\rm
 The $J_\alpha$'s are sometimes called $g$-related almost complex structures (cf. \cite{DLZ1,DLZ2}). %
 Note that the fibre bundle $\mathcal Z$ of all $g$-compatible almost complex structures is called %
 the twistor space of $(M,g,J)$ (see \cite{AHS}):
 $$ \mathcal Z = \{ J \in {\rm SO}(TM) \mid J^2=-{\rm id} \} = {\rm SO}(TM) / {\rm U}(2); $$
 so the twistor fibration $\pi: \mathcal Z \rightarrow M$ is an $S^2$ bundle.
 By using $g$-related almost complex structures one can study
 the twistor space of Riemannian 4-manifolds, complex structures on Riemannian 4-manifolds, %
 Gromov-Witten invariants for K\"{a}hler surfaces and the dimension of
 $J$-anti-invariant cohomology of closed almost complex 4-manifolds
 (cf. \cite{DLZ1,DLZ2,Lee}). } %
 \end{rem}


 \section{Proof of Theorem \ref{thm1}}\setcounter{equation}{0}

 In this section we prove Theorem \ref{thm1}.
 Let us first describe the $C^\infty$-topology
 on the space $\mathcal {J}^\infty$ of $C^\infty$ almost complex structures on $M$.
 For $k=0,1,2,\cdots$, the space $\mathcal{J}^k$ of $C^k$ almost complex structures on $M$
 has a natural separable Banach manifold structure.
 The natural $C^\infty$-topology on $\mathcal{J}^\infty$ is induced by
 the sequence of $C^k$ semi-norms $\|\cdot\|_k$, $k=0,1,2,\cdots$.
 With this $C^\infty$-topology, $\mathcal{J}^\infty$ is a Fr\'{e}chet manifold.
 A complete metric which induces the $C^\infty$-topology on $\mathcal{J}^\infty$
 is defined by
 $$ d(J_1,J_2)=\sum_{k=0}^\infty\frac{\|J_1-J_2\|_k}{2^k(1+\|J_1-J_2\|_k)}. $$
 For details, see \cite{Au,DLZ2}. %


 \vskip 12pt

 \noindent{\bf Proof of Theorem \ref{thm1}.}
 Let $(M,g,J,F)$ be a closed almost Hermitian 4-manifold.
 Note that $H_J^-\subset \mathcal {H}_g^+$ and hence $h_J^-\leq b^+$.
 We assume $b^+ \ge 2$ since the case where $b^+ = 1$
 has been proved by Draghici, Li and Zhang (cf. \cite[Theorem 3.1]{DLZ2}). %

 \vskip 6pt

 To prove the denseness statement,
 we may consider a family $J_t$, $t \in (0,1)$ of almost complex structures on $M$ which is a deformation of $J$, %
 that is, $J_t \rightarrow J$ in the $C^{\infty}$-topology as $t \rightarrow 0$.

 If $h_J^-=0$, then as noted in \cite{DLZ2}, we can establish path-wise semi-continuity property
 for $h_J^{\pm}$
 which follows directly from Lemma \ref{31} and a classical result of Kodaira and Morrow (\cite[Theorem 4.3]{KoMo})
 showing the upper semi-continuity of the kernel of a family of elliptic differential operators.
 Therefore $h_{J_t}^-=0$ for small $t$.

 \vskip 6pt

 We now assume that $h_J^- \ge 1$. Let us write $m: = h_J^-$ and $l:= b^+-m$. %
 We shall construct a family of $g$-compatible almost complex structures
 $\{J_c\} \subset \mathcal {J}_{g}^{\rm c}$ where $c$ are cut-off functions
 to be chosen such that $h_{J_c}^-=0$ and
 $J_c \rightarrow J$ in the $C^{\infty}$-topology as $c \rightarrow 0$.

 First, suppose that $m \,(= h_J^-) < b^+$.
 Then $H_J^{-,\perp}\neq\emptyset$ and
 \begin{align}
 l ={\rm dim}\,H_J^{-,\perp} = b^+-m \ge 1.
 \end{align}
 For each nonzero $[\omega] \in H_J^{-,\perp}$ where $\omega$ is a self-dual harmonic 2-form, we set
 \begin{align}
 f_{\omega} := \langle \omega, F \rangle \in C^{\infty}(M). 
 \end{align}
 Then, by Lemma \ref{32}, $f_{\omega} \not\equiv 0$.
 Set
 \begin{align}
 S_J := \{ \omega \in H_J^{-,\perp} \mid \textstyle\int_M \omega^2 = 1 \}.
 \end{align}
 Then $S_J$ is a sphere of dimension $b^+ -m-1=l-1$.
 Define a function $V:S_J\rightarrow \mathbb{R}$ as follows: for any
 $\omega \in S_J$,
 \begin{align}
 V(\omega) := {\rm vol}\,(M\setminus f_\omega^{-1}(0)) = \textstyle\int_{M\setminus f_{\omega}^{-1}(0)}d\mu_g.
 \end{align}
 Then $V$ is a continuous function on $S_J$. When $H_J^{-,\perp}$ is non-empty,
 \begin{align}\label{deltaJ}
 \delta_J := \inf_{\omega \in S_J} V(\omega) > 0
 \end{align}
 since $S_J$ is compact. Let $\alpha_1,\cdots,\alpha_m \in \mathcal Z_J^-$ be such that
 $[\alpha_1],\cdots,[\alpha_m]$ is an orthonormal basis of $H_J^-$ with respect to the cup product. %
 Choose
 $$ [\omega_1], \cdots, [\omega_l] \in  H_J^{-,\perp}, \quad l = b^+ - m $$
 such that
 $[\alpha_1],\cdots,[\alpha_m],[\omega_1],\cdots,[\omega_l]$
 form an orthonormal basis of
 $\mathcal {H}_g^+$ with respect to the cup product; namely, for $1\leq i \leq m$ and $1\leq j \leq l$, %
 \begin{align}
 \textstyle\int_M\omega_i^2=\int_M\alpha_j^2=1, \quad \int_M\omega_{i}\wedge\alpha_{j}=0, %
 \end{align}
 and, for $1\leq i_1\neq i_2\leq m$ and $1\leq j_1\neq j_2 \leq l$, %
 \begin{align}
 \textstyle\int_M\omega_{i_1}\wedge\omega_{i_2}=0,
 \int_M\alpha_{j_1}\wedge\alpha_{j_2}=0. %
 \end{align}


 To complete the proof of the denseness statement in Theorem \ref{thm1},
 we need the following lemma which is a special case of a theorem of C. B\"{a}r.

 \begin{lem}\label{33}{\rm (\cite[Main Theorem]{Ba})}
 Let $M$ be a closed Riemannian 4-manifold.
 Then every harmonic 2-form $\alpha$ on $M$ has the unique continuation property.
 Hence if $\alpha \not\equiv 0$, then its nodal set $\alpha^{-1}(0)$ has empty interior; %
 in fact, $\alpha^{-1}(0)$ has Hausdorff dimension $\leq 2$.
 \end{lem}

 By Lemma \ref{33}, since $\alpha_1,\cdots,\alpha_m \in \mathcal Z_J^-$,
 the set $\bigcup_{i=1}^{m}\alpha_i^{-1}(0)$ has Hausdorff dimension $\leq2$. %
 Hence $$ M':=\textstyle\bigcap_{i=1}^m(M\setminus \alpha_i^{-1}(0)) $$
 is an open submanifold of $M$ of full volume: ${\rm vol}(M')={\rm vol}(M)$.
 Choose an open set $U\subset M'$ such that ${\rm vol}(U)<\delta_J$.
 Then $\alpha_j|_{U}$, $1\leq i\leq m$ are nonzero sections of $\Lambda_J^-|_{U}$.

 We now construct a new $g$-compatible almost complex structure on $M$ (cf. \cite{DLZ1,Lee}).
 Choose a cut-off function $c_1$ such that $\overline{{\rm supp}\, c_1}\subset U$
 and
 \begin{align}
 |c_1\alpha_1| < 1.
 \end{align}
 Then, by Proposition \ref{prop3.1}, we get %
 \begin{align}
 && J_{c_1}:=J_{c_1\alpha_1}
 = \frac{1-|c_1\alpha_1|^2}{1+|c_1\alpha_1|^2} J - \frac{2}{1+|c_1\alpha_1|^2} K_{c_1\alpha_1}, \\ %
 && F_{c_1}:=F_{c_1\alpha_1}
 = \frac{1-|c_1\alpha_1|^2}{1+|c_1\alpha_1|^2} F + \frac{2}{1+|c_1\alpha_1|^2} c_1\alpha_1. %
 \end{align}
 Thus it is easy to see that, as $c_1 \rightarrow 0$,
 \begin{align}
 J_{c_1} \rightarrow J, \quad F_{c_1} \rightarrow F. %
 \end{align}
 Note that $J_{c_1} \in \mathcal {J}_{g}^{\rm c}\cap \mathcal {J}_{F}^{\rm t}$ %
 and $(g, J_{c_1}, F_{c_1})$ are a family of almost Hermitian structures on $M$. %



 We claim that $H_{J_{c_1}}^{-} \subset H_{J}^{-}$.
 Since $J_{c_1}$ is $g$-compatible, we have
 \begin{align}
 H_{J_{c_1}}^{-} \subset \mathcal H_g^{+} = H_J^{-} \oplus H_J^{-,\perp}.
 \end{align}
 Given any nonzero $\beta \in \mathcal Z_{J_{c_1}}^{-}$,
 there exist real constants $\xi_i$ and $\eta_j$ such that
 \begin{align}
 \beta=\textstyle\sum_{i=1}^{l}\xi_i\omega_i+\sum_{j=1}^{m}\eta_j\alpha_j,
 \end{align}
 where $1\leq i\leq l$ and $1\leq j\leq m$.
 Without loss of generality, we may assume that
 \begin{align}
 \textstyle\int_M\beta^2 = \sum_{i=1}^{l}\xi_i^2 + \sum_{j=1}^{m}\eta_j^2 = 1.
 \end{align}
 Obviously, $\langle\beta,F_{c_1}\rangle=0$.
 Restricted to $M\setminus \overline{{\rm supp}\, {c_1}}$, we have $F_{c_1}=F$.
 On $M\setminus \overline{{\rm supp}\, c_1}$, we get
 \begin{align}\label{4.16}
 \langle\beta,F_{c_1}\rangle |_{M\backslash \overline{{\rm supp}\, c_1}}
 =\textstyle\sum_{i=1}^{l}\xi_i\langle\omega_i,F\rangle |_{M\backslash \overline{{\rm supp}\, c_1}}
 =0.
 \end{align}
 If $\sum_{i=1}^{l}\xi_i\omega_i\in  H_J^{-,\perp}$ is nontrivial,
 then we put
 \begin{align}\label{4.17}
 \beta_1 = \frac{\sum_{i=1}^{l}\xi_i\omega_i}{\big(\int_M(\sum_{i=1}^{l}\xi_i\omega_i)^2\big)^{1/2}}.
 \end{align}
 Obviously, $\beta_1\in S_J$.
 By (\ref{4.16}) and (\ref{4.17}), $f_{\beta_1}=\langle\beta_1,F\rangle \equiv 0$, %
 restricted to $M\setminus \overline{{\rm supp}\, c_1}$.
 Hence $M\setminus f_{\beta_1}^{-1}(0)\subset \overline{{\rm supp}\, c_1}\subset U$.
 It follows that
 \begin{align}
 V(\beta_1) = {\rm vol}(M\setminus f_{\beta_1}^{-1}(0)) \leq {\rm vol}(U) < \delta_J,
 \end{align}
 contradicting the definition of $\delta_J$ (see (\ref{deltaJ})).
 Therefore $\xi_i=0$ for $1\leq i\leq l$ and $\beta=\textstyle\sum_{j=1}^m\eta_j\alpha_j$.
 Thus we have proved that if $H_J^{-} \neq \emptyset$ then $H_{J_{c_1}}^{-} \subset H_J^{-}$.

 \begin{rem}\label{dim}{\rm
 In \cite{DLZ2}, Draghici, Li and Zhang have considered $g$-related almost complex structures: %
 if almost complex structures $J$ and $\tilde{J}$ are $g$-related then
 $\Lambda_{J}^{-} + \Lambda_{\tilde{J}}^{-} \subset \Lambda_{g}^{+}$ and hence
 $H_{J}^{-} + H_{\tilde{J}}^{-} \subset H_{g}^{+}$. }
 \end{rem}

 Secondly, if $H_J^-=\mathcal H^+$ then we construct any $J_{c_1}$ and $F_{c_1}$ such that %
 $H_{J_{c_1}}^- \subset H_J^-$.
 In summary, we have obtained the following

 \begin{prop}\label{prop4.3}
 Let $(M,g,J,F)$ be a closed almost Hermitian 4-manifold.
 If $h_J^-\ge 1$, then we can construct a $g$-compatible almost complex structure $J_{c_1}$ %
 where the volume of ${\rm supp}\,c_1$ is small enough, so that $H_{J_{c_1}}^- \subset H_J^-$. %
 \end{prop}

 The following observation is the key for the computation of $h_{J_{c_1}}^{-}$.

 \begin{prop}\label{prop4.4}{\rm (\cite[Proposition 3.7]{DLZ2})}
 Suppose $J$ and $\tilde{J}$ are $g$-related almost complex
 structures on a closed 4-manifold $M$, with $\tilde{J} \not\equiv \pm J$.
 Then $$ {\rm dim}\,(H_{J}^{-} \cap H_{\tilde{J}}^{-}) \le 1. $$
 \end{prop}

 Let us return to the proof of the denseness statement in Theorem \ref{thm1}. 
 By Propositions \ref{prop4.3} and \ref{prop4.4},
 $h_{J_{c_1}}^{-} = {\rm dim}\,H_{J_{c_1}}^{-} = {\rm dim}\,(H_{J_{c_1}}^{-} \cap H_J^{-}) \le 1$.

 \vskip 6pt

 Without loss of generality, we may suppose that $h_{J_{c_1}}^{-}=1$.
 Choose $[\alpha] \in H_{J_{c_1}}^{-}$ such that
 $\int_{M}\alpha^2=1$.
 We then have ${\rm dim}\,H_{J_{c_1}}^{-,\perp}=b^+-1$, and
 by the same reason as that for (\ref{deltaJ}),
 \begin{align}
 \delta_{J_{c_1}} := {\rm inf}_{\omega\in S_{J_{c_1}}} V(\omega) > 0. %
 \end{align}
 Choose a cut-off function $c_2$ such that $\overline{{\rm supp}\,c_2}\subset M\backslash\alpha^{-1}(0)$
 (by Lemma \ref{33}) and that
 \begin{align}
  {\rm vol}\,({\rm supp}\,c_2) < \delta_{J_{c_1}}.
 \end{align}
 Construct $$ F_{c_2}=f_1 F_{c_1} + c_2 \alpha $$ such that
 $F_{c_2} \wedge F_{c_2} = 2d\mu_g$. 
 It is easy to see that $J$ and $J_{c_2}$ are both $g$-compatible ($g$-related); %
 thus $J_{c_2} \in \mathcal {J}_{g}^{\rm c} \cap \mathcal {J}_{F_{c_1}}^{\rm t}$.

 We claim that $h_{J_{c_2}}^{-}=0$. Otherwise, there exists nonzero $\beta \in H_{J_{c_2}}^{-}$. %
 Then, by the above construction and by Proposition \ref{prop4.3},
 $H_{J_{c_2}}^{-} \subset H_{J_{c_1}}^{-}$;
 hence $\beta\in H_{J_{c_1}}^{-}$, $\beta=\eta\alpha$, $\eta \neq 0$, and
 \begin{align}
 \langle \beta, F_{c_2} \rangle \equiv 0.
 \end{align}
 Restricted to ${\rm supp}\,c_2$, we have $F_{c_2}=f_1F_{c_1}+c_2\alpha$. %
 It follows that
 \begin{align}
 0 = \langle \beta, F_{c_2} \rangle |_{{\rm supp}\,c_2}
   = \langle \eta\alpha, c_2\alpha \rangle = \eta c_2|\alpha|^2. %
 \end{align}
 Thus $\alpha^{-1}(0) \supset {\rm supp}\,c_2$, contradicting Lemma \ref{33}. %
 Hence $h_{J_{c_2}}^-=0$.

 This completes the proof of the denseness statement in Theorem \ref{thm1}.


 \vskip 6pt

 It remains to prove the openness statement in Theorem 1.1. %
 Since the case $b^+=1$ is proved in \cite{DLZ2} by Draghici, Li and Zhang, we consider the case $b^+\geq 2$.
 Suppose that $J_k \rightarrow J$ in the Fr\'{e}chet space as $k \rightarrow \infty$ and that %
 $$ m_k: = h_{J_k}^- \geq 1. $$
 We need to prove that $h_{J}^- \geq 1$.
 Suppose $h_{J}^-=0$. Then $b^+ = h_J^+ - b^-$. %
 Let $g$ be a $J$-compatible metric and set $F=g(J\cdot,\cdot)$.
 Let $\psi^1, \cdots, \psi^{b^+} \in\mathcal {H}_g^+$ be an orthonormal basis with respect to the cup product, that is,
 \begin{align}
 \textstyle\int_M \psi^i \wedge \psi^j = \delta_{ij}, \quad 1 \le i,j \le b^+.
 \end{align}
 Note that
 $\mathcal{H}_g^+ = H_J^- \oplus  H_J^{-,\perp}$.
 Since $h_J^-=0$, by Lemma 2.6,
 \begin{align}\label{418}
 \psi^i = f^i F + P_J^-d\theta^i,
 \end{align}
 where $\theta^i\in\Omega^1$. %
 Set
 $$ g_k=\textstyle\frac{1}{2}(g(\cdot,\cdot)+g(J_k\cdot,J_k\cdot)), \quad F_k=g_k(J_k\cdot,\cdot). $$ %
 Then $(g_k,J_k,F_k) \rightarrow (g,J,F)$ as $k \rightarrow \infty$.
 Since $m_k = h_{J_k}^- \geq 1$, we may choose an orthonormal basis of $\mathcal{H}_{g_k}^+$ with respect to the cup product as follows:
 $$ \{\omega^l(k)\}_{1 \le l \le b^+-m_k} \cup \{\alpha^l(k)\}_{b^+-m_k+1 \le l \le b^+},$$
 where $\{\alpha^l(k)\}_{b^+-m_k+1 \le l \le b^+}\subset H_{J_k}^-$ and $\{\omega^l(k)\}_{1 \le l \le b^+-m_k}\subset H_{J_k}^{-,\perp}$.
 Denote by $\triangle_{g_k}$ the Hodge-de Rham Laplace operator associated to $g_k$
 and by ${\mathbb G}^k$ the Green operator associated to $\triangle_{g_k}$.
 Then, as done in \cite{DLZ2,Le2},
 \begin{align}
 \psi^i=(\psi^i-{\mathbb G}^k(\triangle_{g_k}\psi^i))+{\mathbb G}^k(\triangle_{g_k}\psi^i)=\psi_{{\rm h},k}^i+\psi_{{\rm ex},k}^i,
 \end{align}
 with $\psi_{{\rm h},k}^i:=\psi^i-{\mathbb G}^k(\triangle_{g_k}\psi^i)$ the $g_k$-harmonic part and
 $\psi_{{\rm ex},k}^i:={\mathbb G}^k(\triangle_{g_k}\psi^i)$ the $g_k$-exact part.
 Thus $\psi_{{\rm h},k}^i\rightarrow\psi^i$ and $\psi_{{\rm ex},k}^i\rightarrow 0$ as $k \rightarrow \infty$. %
 Moreover, if $(\psi_{{\rm h},k}^i)^+$ denotes the $g_k$-self-dual part of $\psi_{{\rm h},k}^i$ and $(\psi_{{\rm h},k}^i)^-$ denotes the $g_k$-anti-self-dual part of $\psi_{{\rm h},k}^i$,
 then we still have $(\psi_{{\rm h},k}^i)^+\rightarrow\psi^i$ as $k \rightarrow \infty$ %
 since $\psi^i$ is $g$-self-dual harmonic form and $g_k\rightarrow g$ as $k\rightarrow\infty$.
 Indeed,
 \begin{eqnarray}\label{439} 
  \int_M\psi^i\wedge\ast_{g_k}\psi^i
  &=& \int_M\mid(\psi_{{\rm h},k}^i)^+\mid^2_{g_k}d\mu_{g_k} +\int_M\mid(\psi_{{\rm h},k}^i)^-\mid^2_{g_k}d\mu_{g_k}\nonumber  \\
  &+& \int_M\mid d\theta^i_k\mid^2_{g_k}d\mu_{g_k}
 \end{eqnarray}
 and
 \begin{eqnarray}\label{440} 
   1 &=& \int_M\psi^i\wedge\psi^i \,\ = \,\ \int_M\psi^i\wedge\ast_g\psi^i \nonumber \\
     &=& \int_M\mid(\psi_{{\rm h},k}^i)^+\mid^2_{g_k}d\mu_{g_k}-\int_M\mid(\psi_{{\rm h},k}^i)^-\mid^2_{g_k}d\mu_{g_k},
 \end{eqnarray}
 where $\psi_{{\rm ex},k}^i=d\theta^i_k$.
 Then, by (\ref{439}) and (\ref{440}),
 \begin{eqnarray}
   & & \hspace{-60pt} \int_M\psi^i\wedge\ast_{g_k}\psi^i-\int_M\psi^i\wedge\ast_g\psi^i \nonumber\\
   &=& 2\int_M\mid(\psi_{{\rm h},k}^i)^-\mid^2_{g_k}d\mu_{g_k}+\int_M\mid d\theta^i_k\mid^2_{g_k}d\mu_{g_k}.
 \end{eqnarray}
 Since $\ast_{g_k}\rightarrow\ast_g$ and $g_k\rightarrow g$ as $k\rightarrow \infty$, we obtain that
 $$ 2\int_M\mid(\psi_{{\rm h},k}^i)^-\mid^2_{g_k}d\mu_{g_k}+\int_M\mid d\theta^i_k\mid^2_{g_k}d\mu_{g_k}\rightarrow 0 $$
 and
 $$ 2\int_M\mid(\psi_{{\rm h},k}^i)^-\mid^2_gd\mu_g+\int_M\mid d\theta^i_k\mid^2_gd\mu_g\rightarrow 0 $$
 as $k\rightarrow\infty$.
 Therefore $(\psi_{{\rm h},k}^i)^+\rightarrow \psi^i$, $(\psi_{{\rm h},k}^i)^-\rightarrow 0$ and
 $\psi_{{\rm ex},k}^i\rightarrow 0$ in $L^2(g)$ as $k\rightarrow\infty$.

 Since $(\psi_{{\rm h},k}^i)^+\in\mathcal{H}_{g_k}^+$, there exist $a_l^i(k) \in \mathbb R$, $1\le l\le b^+$ such that
 \begin{align}
 (\psi_{{\rm h},k}^i)^+ = \sum_{l=1}^{b^+-m_k}a_l^i(k)\omega^l(k)
                  + \sum_{l=b^+-m_k+1}^{b^+}a_l^i(k)\alpha^l(k).
 \end{align}
 We have obtained that $(\psi_{{\rm h},k}^i)^+\rightarrow\psi^i$ in $L^2(g)$ as $k\rightarrow\infty$, so
 $$ \int_M\mid(\psi_{{\rm h},k}^i)^+\mid^2_{g_k}d\mu_{g_k}\rightarrow 1 $$
 as $k\rightarrow\infty$. Then $\{(\psi_{{\rm h},k}^i)^+\}$ are bounded in $L^2(g_k)$.

 Let us write
 $$ (\psi_{{\rm h},k}^i)^{+,'}:=\sum_{l=1}^{b^+-m_k}a_l^i(k)\omega^l(k), \quad
    (\psi_{{\rm h},k}^i)^{+,''}:=\sum_{l=b^+-m_k+1}^{b^+}a_l^i(k)\alpha^l(k). $$
 It is easy to see that
 \begin{align*} 
   \int_M\mid(\psi_{{\rm h},k}^i)^+\mid^2_{g_k}d\mu_{g_k}
   = \int_M\mid(\psi_{{\rm h},k}^i)^{+'}\mid^2_{g_k}d\mu_{g_k}
   + \int_M\mid(\psi_{{\rm h},k}^i)^{+,''}\mid^2_{g_k}d\mu_{g_k}.
 \end{align*}
 Since $\{(\psi_{{\rm h},k}^i)^+\}$ are bounded in $L^2(g_k)$, we have
 $\{(\psi_{{\rm h},k}^i)^{+,'}\}$ and $\{(\psi_{{\rm h},k}^i)^{+,''}\}$ are bounded in $L^2(g_k)$.
 As $g_k$-self-dual-harmonic forms, $\{(\psi_{{\rm h},k}^i)^{+,'}\}$ and $\{(\psi_{{\rm h},k}^i)^{+,''}\}$ are bounded in $L^2_2(g_k)$.
 It follows that $\{(\psi_{{\rm h},k}^i)^{+,'}\}$ and $\{(\psi_{{\rm h},k}^i)^{+,''}\}$ are bounded in $L^2_2(g)$
 since $g_k \rightarrow g$ in the $C^\infty$-Fr\'{e}chet space as $k\rightarrow\infty$.
 Hence we can choose a subsequence $\{(\psi_{{\rm h},k_1}^i)^{+,''}\}$ of $\{(\psi_{{\rm h},k}^i)^{+,''}\}$ such that
 $$ \Omega^-_{J_{k_1}} \ni (\psi_{{\rm h},k_1}^i)^{+,''} \rightarrow (\psi^i_{h,\infty})^{+,''}\in \Omega^-_J $$
 in $L^2_1(g)$ as $k_1\rightarrow\infty$ since $J_{k_1} \rightarrow J $ and $\Omega^-_{J_{k_1}} \rightarrow \Omega^-_J$.
 Since $d(\psi_{{\rm h},k_1}^i)^{+,''}=0$, we have $d(\psi^i_{h,\infty})^{+,''}=0$.
 By the assumption that $h^-_J=0$, we can get $(\psi^i_{h,\infty})^{+,''}=0$ by Lemma \ref{31}.
 Thus $(\psi_{{\rm h},k_1}^i)^{+,''}\rightarrow 0$ and $(\psi_{{\rm h},k_1}^i)^{+,'}\rightarrow\psi^i$ in $L^2(g)$ as $k_1\rightarrow\infty$.
 Then
 \begin{align}
 \int_M{({\psi_{{\rm h},k_1}^i})}^{+,'}\wedge{({\psi_{{\rm h},k_1}^i})}^{+,'}\rightarrow\int_M\psi^i\wedge\psi^i=1 %
 \end{align}
 as $k_1\rightarrow\infty$. On the other hand,
 \begin{align}
 {(\psi_{{\rm h},k_1}^i)}^+ = {({\psi_{{\rm h},k_1}^i})}^{+,'}+ \sum_{l=b^+-m_{k_1}+1}^{b^+}a_l^i(k_1)\alpha^l(k_1).
 \end{align}
 It follows that
 \begin{eqnarray*}
 \int_M{({\psi_{{\rm h},k_1}^i})}^+\wedge{({\psi_{{\rm h},k_1}^i})}^+
 &=&\int_M{({\psi_{{\rm h},k_1}^i})}^{+,'}\wedge{({\psi_{{\rm h},k_1}^i})}^{+,'} \\
 &+&\int_M\sum_{l=b^+-m_{k_1}+1}^{b^+}a_l^i(k_1)\alpha^l(k_1)\wedge\sum_{l=b^+-m_{k_1}+1}^{b^+}a_l^i(k_1)\alpha^l(k_1) \\
 &=&\int_M{({\psi_{{\rm h},k_1}^i})}^{+,'}\wedge{({\psi_{{\rm h},k_1}^i})}^{+,'}+\sum_{l=b^+-m_{k_1}+1}^{b^+}(a_l^i(k_1))^2 \\
 &\ge& \int_M{({\psi_{{\rm h},k_1}^i})}^{+,'}\wedge{({\psi_{{\rm h},k_1}^i})}^{+,'}+(a_{b^+}^i(k_1))^2,
 \end{eqnarray*}
 since $h_{J_{k_1}}^- \geq 1$ by the assumption.
 Note that
 \begin{align*}
 \int_M{({\psi_{{\rm h},k_1}^i})}^+\wedge{({\psi_{{\rm h},k_1}^i})}^+\rightarrow 1, \quad%
 \int_M{({\psi_{{\rm h},k_1}^i})}^{+,'}\wedge{({\psi_{{\rm h},k_1}^i})}^{+,'}+\rightarrow 1
 \end{align*}
 as $k_1\rightarrow\infty$.
 Hence
 \begin{align}\label{429}
 (a_{b^+}^i(k_1))^2\rightarrow0
 \end{align}
 as $k_1\rightarrow\infty$.
 Note that
 $\int_M\psi^i\wedge\psi^j=\delta_{ij}$
 and
 ${({\psi_{{\rm h},k_1}^i})}^+\rightarrow \psi^i$ as $k_1\rightarrow\infty$, for $1\leq i \leq b^+$.
 Thus,
 \begin{align}
 \int_M{({\psi_{{\rm h},k_1}^i})}^+\wedge{({\psi_{{\rm h},k_1}^j})}^+
 = \sum_{l=1}^{b^+}a_l^i(k_1)a_l^j(k_1)\rightarrow\delta_{ij}
 \end{align}
 as $k_1\rightarrow\infty$.
 Denote by $A(k_1)$ the $b^+ \times b^+$ matrix: $A(k_1)=(a_l^i(k_1))$.
 Then
 \begin{align}
 B(k_1)=A(k_1)A(k_1)^t=\big(\textstyle\sum_{l=1}^{b^+}a_l^i(k_1)a_l^j(k_1)\big) \rightarrow I_{b^+}
 \end{align}
 as $k_1 \rightarrow +\infty$. Thus $A(k_1)\rightarrow A\in {\rm O}(b^+)$ as $k_1 \rightarrow \infty$.
 Hence $A^t \in {\rm O}(b^+)$ and %
 \begin{align}
 \textstyle\sum_{i=1}^{b^+} a_l^i(k_1)a_m^i(k_1) \rightarrow \delta_{lm} %
 \end{align}
 as $k_1\rightarrow\infty$. In particular, as $k_1\rightarrow\infty$, %
 \begin{align}
 \textstyle\sum_{i=1}^{b^+} a_{b^+}^i(k_1)^2 \rightarrow 1. %
 \end{align}
 This contradicts \eqref{429}.
 Thus we get $h_J^-\geq1$.

 Thus the almost complex structures $J$ with $h_J^-=0$
 form an open subset of $\mathcal J$ in the Fr\'{e}chet space.
 This completes the proof of Theorem 1.1. \qed

 \vskip 12pt

 \noindent{\bf Acknowledgements.}\,
 The authors would like to thank J. Draghici and T.-J. Li for stimulating email discussions
 and T.-J. Li and K.-C. Chen for their help in sending reference papers.
 The second author would like to thank East China Normal University and Qing Zhou
 for hosting his visit in the fall semester in 2011.


 \vskip 24pt

 \noindent Qiang Tan\\
 School of Mathematical Sciences, Yangzhou University, Yangzhou, Jiangsu 225002, China\\
 e-mail: tanqiang1986@yahoo.com.cn\\

 \vskip 6pt

 \noindent Hongyu Wang\\
 School of Mathematical Sciences, Yangzhou University, Yangzhou, Jiangsu 225002, China\\
 e-mail: hywang@yzu.edu.cn\\

 \vskip 6pt

 \noindent Ying Zhang\\
 School of Mathematical Sciences, Soochow University, Suzhou, Jiangsu 215006, China\\
 e-mail: yzhang@suda.edu.cn\\

 \vskip 6pt

 \noindent Peng Zhu\\
 School of Mathematical Sciences, Yangzhou University, Yangzhou, Jiangsu 225002, China\\
 e-mail: zhupeng2004@126.com

 \end{document}